\documentclass[11pt,reqno,tbtags]{amsart}
\usepackage{amssymb}
\usepackage[latin1]{inputenc}

\newtheorem{theorem}{Theorem}[section]
\newtheorem{lemma}[theorem]{Lemma}

\newtheorem{fact}{Fact}

\theoremstyle{definition}
\newtheorem{example}[theorem]{Example}
\newtheorem*{definition}{Definition}
\newtheorem*{problem}{Problem}
\newtheorem{remark}[theorem]{Remark}
\newtheorem*{remark*}{Remark}
\newtheorem*{acks}{Acknowledgements}

\newcommand{\refT}[1]{Theorem~\ref{#1}}

\newcommand{\refL}[1]{Lemma~\ref{#1}}
\newcommand{\refR}[1]{Remark~\ref{#1}}
\newcommand{\refS}[1]{Section~\ref{#1}}

\newcommand{\refE}[1]{Example~\ref{#1}}
\newcommand{\refF}[1]{Fact~\ref{#1}}
\newcommand{\refand}[2]{\ref{#1} and~\ref{#2}}

\newenvironment{romenumerate}{\begin{enumerate}
 }{\end{enumerate}}

\newcounter{thmenumerate}
\newenvironment{thmenumerate}
{\setcounter{thmenumerate}{0}%
 \def\item{\ifnum\thethmenumerate=0\else\par\fi 
 \addtocounter{thmenumerate}{1}\textup{(\roman{thmenumerate})\enspace}}}
{}

\newcommand\ie{\textit{i.e.}\spacefactor=1000}
\newcommand\eg{\textit{e.g.}\spacefactor=1000}
\newcommand\viz{\text{viz.}\spacefactor=1000}

\newcommand{\as}{a.s.\spacefactor=1000}

\newcommand\ga{\alpha}
\newcommand\gb{\beta}

\newcommand\gam{\gamma}

\newcommand\eps{\varepsilon}

\newcommand\bbR{\mathbb R}

\newcommand\E{\operatorname{\mathbb E{}}}
\renewcommand\P{\operatorname{\mathbb P{}}}

\newcommand\ett{\boldsymbol1}
\newcommand\sign{\operatorname{sign}}

\newcommand\set[1]{\ensuremath{\{#1\}}}
\newcommand\bigset[1]{\ensuremath{\bigl\{#1\bigr\}}}
\newcommand\Bigset[1]{\ensuremath{\Bigl\{#1\Bigr\}}}

\newcommand\bigpar[1]{\bigl(#1\bigr)}
\newcommand\Bigpar[1]{\Bigl(#1\Bigr)}

\newcommand\lrpar[1]{\left(#1\right)}

\newcommand\Bigabs[1]{\Bigl|#1\Bigr|}

\def\rompar(#1){\textup(#1\textup)}    
\newcommand\xfrac[2]{#1/#2}
\newcommand\parfrac[2]{\Bigpar{\frac{#1}{#2}}}

\newcommand\half{\frac12}

\newcommand\ff{\check f}
\newcommand\A[1]{{\rm(A#1)}}

\newcommand{\dd}{\,\textup{d}}

\newcommand{\GW}{Galton--Watson}

\newcommand{\cGWt}{conditioned \GW{} tree}
\newcommand{\Holder}{H\"older}
\newcommand\CS{Cauchy--Schwarz}
\newcommand\CSineq{\CS{} inequality}
\newcommand{\maple}{\texttt{Maple}}

\newcommand\urladdrx[1]{{\urladdr{\def~{{\tiny$\sim$}}#1}}}

\begin{document}

\author{James Allen Fill}
\thanks{Research
supported by NSF grants DMS-0104167 and DMS-0406104 and by The Johns Hopkins
University's Acheson~J.~Duncan Fund for the Advancement of Research
in Statistics}

\address{Department of Applied Mathematics and Statistics,
The Johns Hopkins University}
\email{jimfill@jhu.edu} 
\urladdrx{http://www.ams.jhu.edu/~fill/}

\author{Svante Janson}
\address{Department of Mathematics, Uppsala University, PO Box 480,
SE-751~06 Uppsala, Sweden}
\email{svante.janson@math.uu.se}
\urladdrx{http://www.math.uu.se/~svante/}

\title[Tail asymptotics for some limit variables]
{Precise logarithmic asymptotics for the right tails of some limit
  random variables for random trees} 



\date{January 9, 2007}

\begin{abstract}
For certain random variables that arise as limits of functionals
of random finite trees, we obtain precise asymptotics for the logarithm
of the right-hand tail.  Our results are based on the facts 
(i) that the random variables we study can be represented as
functionals of a Brownian excursion and (ii) that a large
deviation principle with good rate function is known explicitly
for Brownian excursion.  Examples include limit distributions
of the total path length and of the Wiener index in conditioned
Galton--Watson trees (also known as simply generated trees).
In the case of Wiener index (where we recover results proved
by Svante Janson and Philippe Chassaing by a different method)
and for some other examples, a key constant is expressed as
the solution to a certain optimization problem, but the
constant's precise value remains unknown.
\end{abstract}

\maketitle

\renewcommand{\=}{:=}
\newcommand{\epsq}{\eps^{1/2}}
\newcommand{\coi}{C[0,1]}
\newcommand{\oi}{[0,1]}
\newcommand{\bbr}{B_{\textup{br}}}
\newcommand{\bex}{B_{\textup{ex}}}
\newcommand{\bexx}{B_{\textup{ex}}^{x}}
\newcommand{\cbm}{C_{\textup{bm}}[0,1]}
\newcommand{\cbmx}{C_{\textup{bm}}}
\newcommand{\cbr}{C_{\textup{br}}[0,1]}
\newcommand{\cbrx}{C_{\textup{br}}}
\newcommand{\cex}{C_{\textup{ex}}[0,1]}
\newcommand{\hbm}{H_{\textup{bm}}}
\newcommand{\hbr}{H_{\textup{br}}}
\newcommand{\hex}{H_{\textup{ex}}}
\newcommand{\kex}{K_{\textup{ex}}}
\newcommand{\ku}{K_{\textup{su}}}
\newcommand{\cboi}{C^\gb[0,1]}
\newcommand{\cbbm}{C^\gb_{\textup{bm}}[0,1]}
\newcommand{\cbex}{C^\gb_{\textup{ex}}[0,1]}

\newcommand{\xT}{T}
\newcommand{\norm}[1]{\|#1\|}
\newcommand{\normoo}[1]{\|#1\|_\infty}
\newcommand{\hx}[1]{\|#1\|_2}
\newcommand{\hh}[1]{\tfrac12\|#1\|_2^2}

\section{Introduction}
\label{Sintro}

Many authors have proved convergence in distribution of various
functionals of various kinds of random trees.
Many have also considered large-deviation estimates and tail bounds.

In this paper, in an attempt to understand several random variables
that arise as such limits, 
we will obtain precise logarithmic asymptotics for their (right-hand) tails.
For example, we will treat the limit distributions of the total
path length and of the Wiener index in conditioned Galton--Watson trees
(a.k.a.\ simply generated trees), where we recover results proved by
Cs\"org\H{o}, Shi and Yor \cite{CsorgoShiYor} and  
Janson and Chassaing \cite{SJ161} by a different method.

The results will be of the ``quasi-Gaussian'' type
$$ \P(X>x)=\exp\left[ - (1 + o(1))x ^2 / (2\gam^2) \right],$$ as $x\to\infty$,
for some positive number $\gamma$ that differs from case to case;
$\gam$ is given as the solution to a variational problem.
In some cases, we can solve the variational problem and find $\gam$
explicitly, while in other cases we only find bounds for $\gam$.
(Note, however, that the asymptotic distributions are not exactly
Gaussian.
Indeed, the examples we study will all be positive random variables.)

Our results are based on 
the fact that the limit random variables we study here
can be represented
as functionals of a (normalized) Brownian excursion; these representations
have been established previously by various authors, and go back to
the theory by Aldous \cite{AldousII,AldousIII} of
the continuum random tree.

\begin{remark}
We consider only tail asymptotics for the limiting random variables
and not for the actual functionals of random trees of a finite
size.  That is, for certain random variables $X_n$ associated with
trees of size $n$ 
with distributional limit $X$,
we find lead-order asymptotics of $\ln \P(X > x)$ as $x$ becomes large, \ie, of
$ \lim_n \ln \P(X_n > x)$.   More interesting would be large deviations for the
sequence $(X_n)$ itself, that is, asymptotics of $\ln \P(X_n > x_n)$
for sequences 
$x_n \to \infty$. 
Such results, however, fall outside the scope of this paper.  Moreover, among
the applications we consider, only in the case of height (our warm-up
\refE{Emax}) 
are such large-deviation results known by any method; see Flajolet et
al.\ \cite{FGOR}.   
\end{remark}

\begin{remark}
Not all limit variables for random trees have quasi-Gaussian tails. 
One well-known example is the total path length in a binary search
tree (under the so-called random permutation model), where Knessl and
Szpankowski~\cite{KnSz}  
give very sharp tail estimates for the limit distribution;
roughly put, they assert that the left tail decays doubly exponentially
and that the right tail decays
exponentially.  But their results rely on several unproven 
regularity assumptions (as noted in their paper), and it is still an
intriguing open problem to verify the
assumptions and prove these results rigorously.

Of course, the family of trees just cited is not simply generated.
A counterexample functional for simply generated families is
the total left path length minus the total right path length in a
uniformly random 
binary tree, 
which is a measure of the asymmetry of the tree.
(Note that the sum is the total path length treated in \refE{Earea} below.)
This difference converges after
suitable scaling to the center of mass of integrated super-Brownian
excursion (ISE), or equivalently 
to the integral of the head of a Brownian snake, 
see \cite{SJ161} or \cite{SJ176}.
For this limit variable
$S$ we have
$ -\ln\bigpar{\P(S>x)} \sim {} \tfrac 34 10^{1/3}\, x^{4/3}$,
as shown in \cite{SJ161}. See also Marckert \cite{Marckert} and 
Bousquet-M\'elou and Janson \cite{SJ185}.
\end{remark}

The main theorem (together with a technical extension) giving large
deviations for functionals of Brownian excursion is stated in
\refS{S:theorem}.  Each application of the main theorem results in a
variational problem; techniques for solving such problems are
discussed in \refS{S:gamma}.  We present applications in \refS{Sex}.
Finally, the main theorem is proved in \refS{Sproof}. 

\section{A general theorem}
\label{S:theorem}

\subsection{Some notation}

We introduce the following notation.
$\coi$ is the usual space of continuous functions on $\oi$ equipped
with the supremum metric 
$\normoo{f-g}:=\sup_t|f(t)-g(t)|$.
We let
\begin{align*}
\cbm&\=\set{f\in \coi:f(0)=0},
\\
\cbr&\=\set{f\in \coi:f(0)=f(1)=0},
\\
\cex&\=\set{f\in \coi:f(0)=f(1)=0,\ f\ge0};
\end{align*}
these are regarded as subsets of $\coi$ equipped with the same metric.
Note that these spaces are closed subspaces of $\coi$ and thus
complete, separable metric spaces.

We further let $B$ be a standard Brownian motion on $[0,1]$,
$\bbr$ a Brownian bridge, and
$\bex$ a standard Brownian excursion; these are random
elements of $\cbm$, $\cbr$, and $\cex$, respectively
(which explains our notation).

Further, let $H$ be the Sobolev space of all absolutely
continuous functions $f\in\coi$ such that 
$\|f'\|_2^2=\int_0^1|f'(t)|^2 \dd t<\infty$.
(The derivative $f'$ exists a.e., and all statements
below about $f'$ for $f\in H$ should be interpreted a.e.)
We define
\begin{align*}
  \hbm&\=H\cap\cbm,
&
  \hbr&\=H\cap\cbr,
&
  \hex&\=H\cap\cex
.
\end{align*}
(The space $\hbm$ is known as the Cameron--Martin space for
Brownian motion, see
\cite{CM} and \cite[Example 8.19]{SJIII}.) 
Similarly, let $K$ be the closed unit ball in $H$, \ie, the set of 
$f\in H$ such that $\|f'\|_2^2=\int_0^1|f'(t)|^2 \dd t\le 1$, and let
\begin{align}\label{k}
  \kex&\=K\cap\cex
.
\end{align}

\subsection{The main result}

We can now state a general theorem for functionals of Brownian
excursions. We give asymptotic results for the distribution itself as
well as for the moment generating function (\ie, for the Laplace transform) and
for the moments. These three results are equivalent (see the proof)
but often useful in 
different situations.

\begin{theorem}\label{T1}
Let $X=\Phi(\bex)$, where $\Phi$ is a continuous nonnegative functional 
on $\cex$ that is positively homogeneous [\ie, $\Phi(tf)=t\Phi(f)$
when $f\in\cex$ and $t\ge0$] and not identically $0$.
Let 
\begin{equation}\label{gamma}
  \gamma:=\max\set{\Phi(f):f\in \kex}.
\end{equation}
Then $0<\gam<\infty$ and 
\begin{alignat}{2}
-  \ln \P(X>x) &\sim \frac{x^2}{2\gamma^2} & &\text{as $x\to\infty$},
\label{a}
\\  
\ln \E e^{tX} &\sim \frac{\gamma^2}{2} t^2 & &\text{as $t\to\infty$},
\label{b}
\\
\bigpar{\E X^r}^{1/r} &\sim \frac{\gamma}{\sqrt e}r^{1/2} & 
\qquad&\text{as $r\to\infty$}.
\label{c}
\end{alignat}
\end{theorem}

\begin{remark}
\label{compact}
$\kex$ is a compact subset of $\cex$, see \refL{Lcomp} below. Hence,
  the maximum in \eqref{gamma} exists and is finite.
\end{remark}

\begin{remark}
It follows from the proof that if the maximum in \eqref{gamma} is
attained at a unique $f_0$, then this $f_0$ is the typical
shape of the $\bex$ giving exceptionally large $X=\Phi(\bex)$, in the
sense
that if $\bexx$ has the conditional distribution of $\bex$ given
$\Phi(\bex)>x$, then $x^{-1}\bexx$ converges in probability,
as $x\to\infty$, to $\Phi(f_0)^{-1} f_0$.

Note that a Brownian excursion a.s.\ does not belong to $H$, since it is
a.s.\ nowhere differentiable. 
Hence $\bexx$ is, for large $x$, with large probability close to
a suitable multiple of $f_0$, but \as{} not exactly equal to it.
\end{remark}

\refT{T1} will be proved in \refS{Sproof}.
We give several applications in \refS{Sex}.

\subsection{An extension}

In one of the  applications in \refS{Sex}, the functional $\Phi$ is
not continuous on $\cex$ and we need an extension (\refT{T1b} below,
also proved 
in \refS{Sproof}) to \Holder{}
spaces. We define, for $0<\gb\le1$, the \Holder{} space $\cboi$ as the
space of all functions $f\in \coi$ such that 
$|f(x)-f(y)|\le C|x-y|^\gb$ for some $C$ and all $x,y\in\oi$;
$\cboi$ is equipped with the metric given by the norm 
$\sup_x|f(x)|+\sup_{x\neq y}|f(x)-f(y)|/|x-y|^\gb$.
Recall that $B,\bbr,\bex\in\cboi$ \as{} for all $\gb<1/2$.
We define $\cbex\=\cboi\cap\cex$, regarded as a subset of
$\cboi$. Note that, for all $\gb\le1/2$, it follows from the 
Cauchy--Schwarz inequality,
as in \eqref{cs} below, that $H\subset\cboi$, and thus $\kex\subset\cbex$.

\begin{theorem}\label{T1b}
If $0<\gb<1/2$, then \refT{T1} remains valid if $\cex$ is replaced by $\cbex$.
\end{theorem}

\section{Finding $\gamma$}\label{S:gamma}

To find $\gamma$ explicitly for the examples in which we have interest, we
begin with some simplifications.
We assume that $\Phi$ is a continuous functional on $\kex$
as in \refT{T1}.
We begin by listing some properties $\Phi$ may have.

\begin{itemize}
\item[\A1] 
$\Phi$ is \emph{symmetric}: if $\ff(x)\=f(1-x)$ then $\Phi(\ff)= \Phi(f)$.
\item[\A2] 
$\Phi$ is \emph{concave}:
  $\Phi\bigpar{\frac12(f+g)}\ge\frac12\Phi(f)+\frac12\Phi(g)$.
(For positively homogeneous $\Phi$ this is equivalent to superadditivity.)
\item[\A3] 
$\Phi$ is \emph{monotone}: $f_1\le f_2$ implies $\Phi(f_1)\le \Phi(f_2)$.
\end{itemize}

Our first of two lemmas shows that if~$\Phi$ has certain of these
properties, then 
the search space for~$f$ maximizing $\Phi(f)$ may be suitably narrowed
from $\kex$.

\begin{lemma}\label{L1}
Let $\Phi$ be a continuous functional defined on $\kex$.
\begin{romenumerate}
\item
If  $\Phi$ is symmetric and concave, then $\max_{\kex} \Phi(f)$ is attained
by an~$f$ which is symmetric ($\ff=f$).
\item
If  $\Phi$ is monotone, then $\max_{\kex} \Phi(f)$ is attained by a
unimodal $f$, \ie, an $f$ such that
$f'\ge0$ on $(0,a)$ and $f'\le0$ on $(a,1)$ for some $a\in(0,1)$.
\item
If  $\Phi$ is symmetric, concave, and monotone, 
then $\max_{\kex} \Phi(f)$ is attained
by a
symmetric unimodal $f$, \ie, an $f$ such that
$f'\ge0$ on $(0,1/2)$ and $f'(x)=-f'(1-x)$ on $(1/2,1)$.
\end{romenumerate}
\end{lemma}

\begin{proof}
(i):
Let $f\in \kex$ maximize $\Phi(f)$, and let $g=\tfrac12(f+\ff)$.
Then $g\in \kex$ is symmetric and, by the assumptions,
$\Phi(g)\ge \tfrac12\Phi(f)+ \tfrac12\Phi(\ff)=\Phi(f)$.
Hence $g$, too, maximizes $\Phi$.

(ii):
Let $f\in \kex$ maximize $\Phi(f)$, and 
define $g\in H$ by $g(0)=0$ and
$$
g'(x)=
\begin{cases} 
|f'(x)| & \mbox{if\ }x<a \\ 
-|f'(x)| & \mbox{if\ }x>a 
\end{cases}
$$
where $a$ is such that $\int_0^a|f'(x)|\dd x = \tfrac12
\int_0^1|f'(x)|\dd x$.
Then $g(1)=g(0)=0$, $g\ge0$, and $\hx{g'}=\hx{f'}$, so $g\in \kex$.
For $x\le a$ we have 
$$
f(x)\le\int_0^x|f'(y)|\dd y = \int_0^xg'(y)\dd y =g(x)
$$
and for $x\ge a$ 
\begin{equation*}
f(x)\le \int_x^1|f'(y)|\dd y = \int_x^1-g'(y)\dd y =g(x),  
\end{equation*}
so $f\le g$ and thus $\Phi(g)\ge\Phi(f)$.
Hence $g$ too maximizes $\Phi$.

(iii): Argue first as for (i) and then as for (ii).
\end{proof}

Our second lemma concerns maximization of $\Phi(f)$ over a 
certain smaller class of functions~$f$ than $\kex$. 

\begin{lemma}\label{L2}
\begin{thmenumerate}
\item
Let $\ku$ be the subset of $\kex$ consisting of symmetric unimodal
functions.
Suppose that $\Phi$ is a continuous 
functional on $\ku$ such that for
some nonnegative function $h\in L^2[0,1/2]$ 
\begin{equation}\label{sofie}
  \Phi(f)=\int_0^{1/2} f'(t) h(t)\dd t,
\qquad f\in \ku.
\end{equation}
Then 
\begin{equation}\label{emma}
\max_{\ku} \Phi = \frac{1}{\sqrt2}\|h\|_{L^2[0,1/2]}
=\left(\frac12\int_0^{1/2} h(t)^2\dd t\right)^{1/2}.
\end{equation}
A maximizing $f\in \ku$ is given by $f'=(2\max_{\ku}\Phi)^{-1}h$ on $[0,1/2]$.
\item	
Suppose that $\Phi$ is a continuous symmetric, concave, monotone
functional on $\kex$ such that \eqref{sofie} holds for
some nonnegative function $h\in L^2[0,1/2]$ and all $f\in \ku$.
Then 
$\max_{\kex} \Phi = \max_{\ku} \Phi$ is given by \eqref{emma}.
\end{thmenumerate}
\end{lemma}
\begin{proof}
(i):
This is immediate by Hilbert space theory, since 
\begin{equation*}
\bigset{f'|_{[0,1/2]} : f\in \ku}
= \Bigset{g\in L^2[0,1/2]: g\ge0 \text{ and } \hx{g}=\xfrac1{\sqrt2}}.
\end{equation*}

(ii): This follows from (i) and \refL{L1}(iii).
\end{proof}

\section{Applications}\label{Sex}

We give several applications of the general theorem to functionals of
interest for random trees. In all cases, the random trees that we
consider are conditioned Galton--Watson trees, also known as simply
generated trees. As is well-known, this includes several important
types of random trees, for example
random planar trees, random labelled trees, and random binary trees
(in each case uniformly distributed over all trees of the given type with a
given number $n$ of vertices). As is shown in the references given below, 
the functionals we study have limit distributions as the size
$n$ of the random trees tends to infinity, after proper normalization.
Moreover, these limit distributions do not depend on the particular
class of random trees (within the class of conditioned Galton--Watson trees)
except for a simple scale factor. In the results below, we therefore will
not usually discuss the random trees.

Moreover, since the asymptotic results always are 
given by \eqref{a}, \eqref{b}, \eqref{c},
we will only give the value of $\gamma$.

\begin{example}[Height and width]\label{Emax}
For 
both the height and the width of a \cGWt, the limit distribution
(after suitable rescaling) is given by the same random variable,
\viz{} $\max\bex$, see Aldous \cite{AldousII} and 
Chassaing, Marckert and Yor
\cite{ChMY}; see also \cite[Section 7]{SJ167}.
The distribution of this random variable is well-known
\cite{Chung,Kennedy}, see \cite{BianePY} for much more information;
in particular,
\begin{equation}\label{max}
\P(\max_t \bex(t) \le x)
= 1+2\sum_{k=1}^\infty(1-4k^2x^2)\exp(-2k^2x^2),
\qquad x>0.
\end{equation}
Hence the asymptotics we obtain from \refT{T1} do not yield anything
new, but they 
serve as a simple warm-up exemplifying our results.

Thus, let $\Phi(f):=\max f$. 
This functional is symmetric and monotone, but not concave.
By \refL{L1}(ii), the maximum is attained for a unimodal $f$, but we
cannot use \refL{L2}.
We can in this case easily argue directly. Let $f \in \kex$.  By the
Cauchy--Schwarz 
inequality,
for every $x\in[0,1]$,
\begin{equation*}
  \begin{split}
f(x)=\tfrac12\int_0^x f'(t)\dd t	
-\tfrac12\int_x^1 f'(t)\dd t	
=\int_0^1 f'(t)\tfrac12\sign(x-t)\dd t
\le \tfrac12\hx{f'}\le\tfrac12,
  \end{split}
\end{equation*}
with equality if $x=1/2$ and $f'(t)=\sign(\half-t)$, \ie,
if $f(x)=x$ for $0\le x\le 1/2$ and 
$f(x)=1-x$ for $1/2\le x\le 1$.
Thus $\gamma=\max_{\kex}\Phi=1/2$, so 
$\ln\P(\max_t \bex(t) > x) \sim -2x^2$, in accordance with \eqref{max}.
\end{example}

\begin{remark*}
(a) The maximizing $f$ in \refE{Emax} is easily seen to be unique.
We guess that the same is true in all examples below, but we have not
checked this.

(b) \refE{Emax} shows that the maximum $\gam$ may be attained on $\ku$
even if $\Phi$ is not concave.

(c) It follows from Theorem~1.2 in Flajolet et al.\ \cite{FGOR} that
the height $H_n$ of  
a conditioned critical Galton--Watson tree with offspring distribution having
variance $\sigma^2$, when normalized to
$X_n := \sigma H_n / (2 \sqrt{n})$, satisfies, for any $c > 0$, the
``zone of convergence'' result 
$$
\P(X_n > x) \sim \P(\max_t \bex(t) > x)\mbox{\ \ uniformly for $x < c
  \sqrt{\log n}$} ; 
$$
and hence that
$$
-\ln \P(X_n > x_n) \sim 2 x^2_n\mbox{\ \ as $n \to \infty$}
$$
provided $x_n \to \infty$ and $x_n = O(\sqrt{\log n})$.  Presumably,
similar such results hold 
for other functionals treated below, but our techniques cannot yield
these more delicate results. 
\end{remark*}

\begin{example}[Total path length] \label{Earea}
It is also well known that the asymptotic distribution of the total
path length in a \cGWt{} is given by the
Brownian excursion area $\int_0^1\bex(t)\dd t$ \cite{AldousII,AldousIII}.
Thus, we
now  let $\Phi(f)=\int_0^1 f$.
Here $\Phi$ is symmetric, concave (in fact, linear), and monotone.
If $f\in \ku$, then by integration by parts,
\begin{equation*}
  \Phi(f)=2\int_0^{1/2} f(t)\dd t
=2\int_0^{1/2} (\tfrac12-t)f'(t)\dd t.
\end{equation*}
Hence \refL{L2} applies with $h(t)=1-2t$, which gives
\begin{equation*}
\gamma=\max_{\ku} \Phi 
=\left(\frac12\int_0^{1/2} (1-2t)^2\dd t\right)^{1/2}
=\frac1{\sqrt{12}}.
\end{equation*}
This agrees with the tail asymptotics given by 
Cs\"org\H{o}, Shi and Yor, 
\cite[Proof of Theorem 3.1]{CsorgoShiYor}.
A maximizing $f$ is given by $f'=h/(2\gamma)=\sqrt3(1-2t)$ on $[0,1/2]$,
and thus $f(t)=\sqrt3t(1-t)$, $t\in\oi$.

The variable $\xi=2\int_0^1 \bex$ is studied in \cite{SJ161}.
\refT{T1} applies with 
$\Phi(f)=2\int_0^1 f$. This is simply twice the functional just studied,
and hence $\gamma =2/\sqrt{12}=1/\sqrt3$, in agreement with the result in
\cite[Remark 4.9]{SJ161} (obtained there by a different method).
A maximizing $f$ is again $f(t)=\sqrt3t(1-t)$.
\end{example}

In the following examples we use the notation
\begin{equation*}
  m(f;s,t)\=\inf\set{ f(u):u\in[s,t]}
\end{equation*}
for a function $f$ on $\oi$ and $0\le s\le t\le 1$.

\begin{example} \label{Eeta}
Another random variable studied
in \cite{SJ161} is
$\eta\=4\iint_{s<t} m(\bex; s,t) \dd s\dd t$;
this arises as the limit for the sum, over all pairs of vertices in
the random tree, of the depth of the last common ancestor. 
\refT{T1} applies with 
$\Phi(f)=4\iint_{s<t} m(f; s,t) \dd s\dd t$.
This $\Phi$ is symmetric, concave, and monotone.
For $f\in \ku$, $m(f;s,t)=f(s)$ when $|\half-s|>|\half-t|$, and
$m(f;s,t)=f(t)$ otherwise. Hence, using the symmetry of $f$,
\begin{equation*}
\begin{split}
\Phi(f)
&=8\hskip-1em\iint\limits_{\substack{s<t\\|\half-s|>|\half-t|}} 
\hskip-1em f(s) \dd s\dd t
=8\int_0^{1/2}\int_s^{1-s}f(s)\dd t\dd s \\
&=8\int_0^{1/2}(1-2s)f(s)\dd s
=2\int_0^{1/2}(1-2s)^2f'(s)\dd s.
\end{split}
\end{equation*}
Thus \refL{L2} applies with $h(t)=2(1-2t)^2$ and
\begin{equation*}
\gamma=\max_{\ku} \Phi 
=\left(\frac12\int_0^{1/2}4 (1-2t)^4\dd t\right)^{1/2}
=\frac1{\sqrt{5}}.
\end{equation*}

This gives a new proof of the result in \cite[Theorem 4.6]{SJ161}.
A maximizing function is given by $f'(t)=h/2\gamma=\sqrt5(1-2t)^2$,
$t\le1/2$,
and thus $f(t)=\frac{\sqrt5}{6}\bigpar{1- |1-2t|^3}$, $t\in\oi$.
\end{example}

\begin{example}[Wiener index] \label{Ezeta}
It is shown in \cite{SJ146} that for the Wiener index of the random tree, the
limit random variable is
$\zeta\=\xi-\eta$, with $\xi$ and $\eta$ given in the preceding examples.
Thus \refT{T1} applies to $\zeta$ with 
$\Phi(f)=2\iint_{s<t} \bigpar{f(s)+f(t)-2m(f; s,t)} \dd s\dd t$.
This $\Phi$ is symmetric, but neither concave (on the contrary, it is
convex) nor monotone. 
For $f\in \ku$, Examples \refand{Earea}{Eeta} show that
\eqref{sofie} holds with $h(t)=2(1-2t)-2(1-2t)^2=4t(1-2t)$.
Hence \refL{L2} shows that
\begin{equation*}
\max_{\ku} \Phi 
=\left(\frac12\int_0^{1/2}16 t^2 (1-2t)^2\dd t\right)^{1/2}
=\frac1{\sqrt{30}}.
\end{equation*}
However, we do not know whether this also is the maximum $\gam$ over
$\kex$,
so we can only conclude $\gam\ge1/\sqrt{30}$.

An upper bound can be found as follows.
If $f\in\kex$,
let $0<s<t<1$ and let $v$ be a minimum point for $f$ in $[s,t]$,
\ie, a point $v\in[s,t]$  such that $f(v)=\min\set{f(u):u\in[s,t]}=m(f;s,t)$.
Then
\begin{align*}
f(s) + f(t) - 2 m(f; s, t)
={}-\int_s^v f'(u)\dd u + \int_v^t f'(u)\dd u
\le \int_s^t |f'(u)|\dd u  .
\end{align*}
Thus, by the \CSineq{} and the assumption $f\in K$,
\begin{align*}
\Phi(f)
&=
2\iint_{s<t}\!\left[ f(s) + f(t) - 2 m(f; s, t) \right]\dd s\dd t \\
 &\leq 2 \iint_{s<t}\!\left[ \int^t_s | f'(u) |\dd u \right]\dd s\dd t \\
 &  =  2\int^1_0\!| f'(u) | \left[ \iint_{0 < s < u < t < 1}\!\dd s\dd t
          \right]\dd u \\
 &  =  2\int^1_0\!| f'(u) |\,u (1 - u)\dd u \\
 &\leq 2\| f' \|_2 \left[ \int^1_0\!u^2 (1 - u)^2\dd u \right]^{1/2} \\
 &  =  2\| f' \|_2 / \sqrt{30} \leq 2 / \sqrt{30}.
\end{align*}
Consequently, $\gam\le2/\sqrt{30}$, and combining this with the lower
bound above we find
$1/\sqrt{30}\le\gam\le2/\sqrt{30}$.

\begin{problem}
  Find $\gam$ for the random variable $\zeta$.
\end{problem}

\end{example}

Fill and Kapur \cite{FK} \cite{FKsg} have studied the sum, over all
vertices~$v$ 
in the random tree, of the $\ga$th power of the size of the subtree
rooted at~$v$; 
here $0 < \ga < \infty$ is a parameter.  For $\ga > 1/2$, which is the
only range 
we shall consider here, they show that,
after suitable scaling, there is a limit distribution characterized by
its moments. 
Let $Y_{\ga}$ have this distribution.  Fill and Janson \cite{FJbexrep}
show that 
$Y_{\ga}$ can be represented as $\Phi(\bex)$ with
\begin{multline}\label{wa}
\Phi(f)=\ga \int^1_0 \left[ t^{\ga - 1} + (1 - t)^{\ga - 1} \right]f(t)\dd t 
\\
 {} - \ga (\ga - 1) \iint\limits_{0 < s < t < 1} 
 (t - s)^{\ga - 2} \left[ f(s) + f(t) - 2 m(f;s,t) \right]\dd s\dd t. 
\end{multline}
Note that for $\ga=1$ this reduces to $2\int_0^1 f$, and thus $W_1=\xi$
in \refE{Earea}. Moreover, if $\ga>1$, then \eqref{wa} simplifies to
\begin{equation}\label{wa1}
\Phi(f)=2\ga (\ga - 1) \iint\limits_{0 < s < t < 1} 
 (t - s)^{\ga - 2}  m(f;s,t) \dd s\dd t. 
\end{equation}
In particular, $W_2=\eta$ in \refE{Eeta}.

\begin{example} \label{Ealpha}
If $\ga>1$, \refT{T1} applies with $\Phi$ given by \eqref{wa1}.
This $\Phi$ is symmetric, concave, and monotone.
For $f\in \ku$, as in \refE{Eeta},
\begin{equation*}
\begin{split}
\Phi(f)
&=4\ga(\ga-1)\hskip-1em\iint\limits_{\substack{s<t\\|1/2-s|>|1/2-t|}} 
\hskip-1em (t-s)^{\ga-2}f(s) \dd s\dd t
\\
&=4\ga(\ga-1)\int_0^{1/2}\int_s^{1-s}(t-s)^{\ga-2} f(s)\dd t\dd s 
\\
&=4\ga\int_0^{1/2}(1-2s)^{\ga-1}f(s)\dd s
=2\int_0^{1/2}(1-2s)^\ga f'(s)\dd s.
\end{split}
\end{equation*}
Thus \refL{L2} applies with $h(t)=2(1-2t)^\ga$ and
\begin{equation*}
\gamma\=\max_{\kex} \Phi 
=\left(\frac12\int_0^{1/2}4 (1-2t)^{2\ga}\dd t\right)^{1/2}
=\frac1{\sqrt{2\ga+1}}.
\end{equation*}
This has been found (in the form \eqref{c}) by Fill and Kapur \cite{FKasymom}.
A maximizing function is given by 
$f'(t)=h(t)/(2\gamma)=\sqrt{2\ga+1}(1-2t)^\ga$,
$t\le1/2$,
and thus 
$f(t)=\frac{\sqrt{2\ga+1}}{2(\ga+1)}\bigpar{1- |1-2t|^{\ga+1}}$,
$t\in\oi$.
\end{example}

\begin{example}
  Now let $1/2<\ga<1$.
In this case, the formula \eqref{wa1} cannot be used (the integral
diverges unless $f$ is constant; moreover, the factor in front is
negative), so we have to use \eqref{wa}.
When $\ga<1$, this functional $\Phi$ is \emph{not} continuous on
$\cex$.
It is, however, continuous on the \Holder{} space $\cbex$ when
$\ga+\gb>1$, as is easily verified.
We thus choose $\gb\in(1-\ga,1/2)$ and use \refT{T1b}.

Nevertheless, there are further problems. When $\ga<1$, the functional
$\Phi$ is 
neither monotone nor concave (it is instead convex), so we cannot
apply \refL{L2}.

For $f\in \ku$ we find, in similar fashion as for \refE{Ealpha},
omitting the details, 
\begin{equation*}
\Phi(f)=2\int_0^{1/2}(1-2s)^\ga f'(s)\dd s,
\end{equation*}
and thus, also for $\ga<1$, 
\begin{equation*}
\max_{\ku} \Phi 
=\frac1{\sqrt{2\ga+1}}.
\end{equation*}
However, for $\ga<1$, we do not know whether this also is the maximum
over $\kex$, so we can only conclude
$\gam\ge1/\sqrt{2\ga+1}$. 

To get an upper bound, assume $f\in\kex$ and denote the two integrals
in \eqref{wa} by 
$\Phi_1(f)$ and $\Phi_2(f)$. An integration by parts yields
\begin{align*}
\Phi_1(f)
 & = \int^1_0\!f'(u) \left[ (1 - u)^{\ga} - u^{\ga} \right]\,du
\le \int^1_0\!|f'(u)| \left| (1 - u)^{\ga} - u^{\ga} \right|\,du,
 \end{align*}
while an argument as in \refE{Ezeta} yields
\begin{align*}
\Phi_2(f)
 &\leq \int^1_0\!| f'(u) | \left[ u^{\ga} + (1 - u)^{\ga} - 1 \right]\,du.
\end{align*}
Hence, if we define
\begin{equation*}
  h(u):=
| (1 - u)^{\ga} - u^{\ga}|
+ u^{\ga} + (1 - u)^{\ga} - 1
=
  \begin{cases}
	2 (1 - u)^{\ga} - 1,& \mbox{if\ }0\le u\le 1/2,\\
	2 u^{\ga} - 1,& \mbox{if\ }1/2\le u\le 1,
  \end{cases}
\end{equation*}
we have, for $f\in\kex$,
\begin{align*}
\Phi(f)
  &  =  \Phi_1(f) + \Phi_2(f) 
  \leq \int^1_0\!| f'(u) | h(u)\,du 
  \leq \| f' \|_2 \left[ \int^1_0\! h(u)^2\,du  \right]^{1/2},
\end{align*}
and thus
\begin{align*}
\gam 
&\le
\lrpar{ \int^1_0\! h(u)^2\,du }^{1/2}
=\lrpar{\frac{8}{2 \ga + 1}(1-2^{-2\ga-1}) 
- \frac{8}{\ga + 1}(1-2^{-\ga-1}) + 1}^{1/2}.
\end{align*}
Denoting the right hand side by $\psi(\ga)^{1/2}$, we have verified (first
graphically using \maple, and then rigorously using calculus)
that $(2\ga+1)\psi(\ga)$ is decreasing on $[1/2,1]$, and thus 
the maximum is attained for $\ga=1/2$, which gives the value
$8(\sqrt2-1)/3$. Hence, for $1/2<\ga<1$,
\begin{equation*}
\gam
\le 
\psi(\ga)^{1/2}
\le \left(\frac{8(\sqrt2-1)}{3(2\ga+1)}\right)^{1/2}
\le \frac{1.051}{\sqrt{2\ga+1}}.  
\end{equation*}
Hence our upper and lower bound differ by a factor less than $1.051$
(and the ratio tends to 1 as $\ga\to1$).

\begin{problem}
  Find $\gam$ for $W_\ga$ when $\ga<1$.
\end{problem}
\end{example}

\section{Proof of Theorems \refand{T1}{T1b}}\label{Sproof}

\begin{proof}[Proof of \refT{T1}]
We begin with a simple lemma, see \eg{} \cite[Lemma 27.7]{Kallenberg}.
\begin{lemma}\label{Lcomp}
  The set $\kex$ defined at~\eqref{k} is a compact subset of $\cex$.
\end{lemma}
\begin{proof}
If $f\in K$ and $0\le x\le y\le 1$, then the Cauchy--Schwarz
inequality yields
  \begin{equation}\label{cs}
|f(x)-f(y)|^2 
=
\Bigabs{\int_x^y f'(t) \dd t }^2
\le \int_x^y \dd t \int_x^y |f'(t)|^2 \dd t 
\le y-x.
  \end{equation}
Since further $f\in\kex$ implies $f(0)$=0, it follows from the
Arzel\`a--Ascoli theorem that $\kex$ is relatively compact in $\coi$,
and thus in $\cex$.

It remains to show that $\kex$ is a closed subset of $\cex$. Thus, assume that
$f_n\in\kex$ and that $f_n\to f$ in $\coi$. The functions $f_n'$
belong to the unit ball of $L^2\oi$, so by weak compactness there
exists a subsequence $f'_{n_k}$ that converges weakly in $L^2$, say to $g$.
Define $F(x)\=\int_0^x g(t)\dd t$. Then $F'=g$ a.e., and $F\in \kex$.
Moreover, 
the weak convergence $f'_{n_k}\to g$ along the subsequence implies,
for every $x\in\oi$, 
\begin{equation*}
  f_{n_k}(x)
=\int_0^x   f'_{n_k}(t) \dd t
=\langle f_{n_k},\ett_{[0,x]} \rangle
\to 
\langle g,\ett_{[0,x]} \rangle
=
\int_0^x   g(t) \dd t
=F(x).
\end{equation*}
Hence $f = F \in\kex$.
\end{proof}

As noted at \refR{compact}, \refL{Lcomp} shows that the maximum $\gam$
in \eqref{gamma} 
exists and is finite. Moreover, $\gam>0$, because otherwise
$\Phi(f)=0$ for every $f\in\kex$. By homogeneity, this would imply
$\Phi(f)=0$ for every $f\in\hex$. However, $\hex$ is dense in $\cex$,
as can be seen by approximating a continuous function by
piecewise linear functions, and since $\Phi$ is assumed to be
continuous, this would imply that $\Phi$ vanishes identically on
$\cex$, contrary to our assumption.

To prove \refT{T1}, 
we use some notations and results from large deviation theory, see
for example Kallenberg \cite[Chapter 27]{Kallenberg}
or 
Dembo and  Zeitouni \cite{DemboZ}.

\begin{definition}[{\cite[pp.\ 545--546]{Kallenberg}}]
  A family $(X_\eps)_{\eps>0}$ of random elements in some metric space
  $S$ satisfies the \emph{Large Deviation Principle (LDP) with good
  rate function $I$}
  if $I:S\to[0,\infty]$ is a function such that the level sets
\set{x\in S:I(x)\le r} are compact for all finite $r$ and,
for every Borel set $A\subseteq S$,
\begin{equation*}
-\inf_{x\in A^\circ} I(x) 
\le
\liminf_{\eps\to0}\bigpar{ \eps \ln \P(X_\eps\in A)}
\le
\limsup_{\eps\to0} \bigpar{\eps \ln \P(X_\eps\in A)}
\le
-\inf_{x\in \overline A} I(x)
. 
\end{equation*}
\end{definition}

We begin with two central facts.

\begin{fact}[{\cite[Theorem 27.6]{Kallenberg}}]
\label{F1}
If $B$ is a Brownian motion, then
$(\eps^{1/2}B)$ satisfies the LDP in $\cbm$ with good rate function
$I(f)=\hh{f'}$ 
for $f\in \hbm$ and $I(f)=\infty$ otherwise.
\end{fact}

\begin{fact}[{\cite[Theorem 27.11]{Kallenberg}}]
\label{F2}
If $F:S\to T$ is a continuous mapping of one metric space into
another, and $X_\eps$ satisfies the LDP in $S$ with good rate function
$I$,
then $F(X_\eps)$ satisfies the LDP in $T$ with good rate function 
$J(y) \= \inf\set{I(x):F(x)=y}$.
\end{fact}

For the first application of \refF{F2}, note that a Brownian bridge may be
constructed by $\bbr(t)\=B(t)-tB(1)$. 
Hence, let $F(f)(t):=f(t)-tf(1)$. This is a
continuous map $\cbm\to\cbr$ and $F(B)= \bbr$.
It is easily seen that 
$J(f)=\inf_{a\in \bbR}\tfrac12\|f'+a\|_2^2=I(f)$ for $f\in\hbr$ 
and $J(f)=\infty$ otherwise.
Hence Facts \refand{F1}{F2} yield the LDP for the Brownian bridge:

\begin{fact}[{\cite[Exercise 27.10]{Kallenberg}}]
\label{F3}
If $\bbr$ is a Brownian bridge, then
$(\eps^{1/2}\bbr)$ satisfies the LDP in $\cbr$ with good rate function
$I(f)=\tfrac12\|f'\|_2^2$ for $f\in \hbr$ and $I(f)=\infty$ otherwise.
\end{fact}

Moreover
\cite{Kallenberg},
Facts \refand{F1}{F3} extend readily
to $d$-dimensional Brownian motion and
bridge, respectively, if we 
replace the spaces $\cbm$, $\hbm$, $\cbr$, and $\hbr$ by
the corresponding spaces
$\cbmx(\oi,\bbR^d)$, and so on,
of functions with values in $\bbR^d$, interpreting 
$\|f'\|_2^2=\int_0^1|f'(t)|^2\dd t$ with $|f'(t)|$ the usual Euclidean
length of the vector $f'(t)$ in $\bbR^d$.

Turning to the Brownian excursion, we use the result 
that $\bex$ has the same distribution as the process $|\bbr^{(3)}|$,
where $\bbr^{(3)}$  
is 3-dimensional
Brownian bridge, see, \eg{}, Revuz and Yor \cite[Theorem XII.(4.2)]{RY}.
We can thus apply \refF{F2} with 
$S=\cbrx(\oi,\bbR^3)$, $T=\cex$, $F(g)=|g|$ and 
$X_\eps=\eps^{1/2} \bbr^{(3)}$. Recalling Fact~\ref{F3} and noting
that $|F(g)'|\le |g'|$, 
it is easily seen that 
$J(f)\=\inf\set{ I(g):g\in \cbrx(\oi,\bbR^3) \text{ and } |g|=f}$
equals $\tfrac12\|f'\|_2^2$ for $f\in \hex$ and equals $\infty$ otherwise, and
we obtain the following result.

\begin{fact}[Serlet \cite{Serlet}]
\label{F4}
If $\bex$ is a standard Brownian excursion, then
$(\eps^{1/2}\bex)$ satisfies the LDP in $\cex$ with good rate function
$I(f)=\tfrac12\|f'\|_2^2$ for $f\in \hex$ and $I(f)=\infty$ otherwise.
\end{fact}

(It is also possible, but more complicated, to prove this from \refF{F3}
using the result by Vervaat \cite{Vervaat} 
that the random process
$\bbr(t)-\min \bbr$ has the same distribution as $\bex(U+t)$, where
$U$ is uniform on $[0,1]$ and independent of $\bex$, and addition is
modulo 1.)

Finally, we apply \refF{F2} once more, now to $\Phi:\cex\to \bbR$
and find that $\epsq X=\Phi(\epsq \bex)$ satisfies the LDP in
$[0,\infty)$ with the good
rate function, for $x>0$, 
\begin{equation*}
  \begin{split}
\inf_{f\in \hex:\; \Phi(f)=x} \hh{f'}
&=   
\inf_{f\in \hex:\; \Phi(f)\neq0} 
  \tfrac12\Bigl\|\frac{xf'}{\Phi(f)}\Bigr\|_2^2\\
&=\inf_{f\in \hex:\; \hx{f'}=1,\,\Phi(f)\neq0} \tfrac12\parfrac{x}{\Phi(f)}^2
=\frac{1}{2\gamma^2}x^2.	
  \end{split}
\end{equation*}
Taking $A=(1,\infty)$ and $\eps=x^{-2}$ in the definition of LDP, 
this proves \eqref{a}.
Finally, \eqref{b} and \eqref{c} follow easily from \eqref{a} by
integration; indeed, the (more difficult) converses hold too, see
Davies \cite{Davies} and Kasahara \cite{Kasahara} or
\cite[Theorem 4.5]{SJ161}.
This completes the proof of \refT{T1}.
\end{proof}

\begin{proof}[Proof of \refT{T1b}]
We begin by observing that the following extension of \refF{F1}
holds, also in $d$ dimensions.
\begin{fact}
\label{F5}
If $0<\gb<1/2$, then
$(\eps^{1/2}B)$ satisfies the LDP in $\cbbm\=\cbm\cap\cboi$ with good
rate function 
$I(f)=\hh{f'}$ 
for $f\in H\cap\cbbm$ and $I(f)=\infty$ otherwise.
\end{fact}
Indeed, by \cite[Theorem 27.11(ii)]{Kallenberg}], this follows from
\refF{F1} and the property that 
$(\eps^{1/2}B)$ is \emph{exponentially tight} in $\cbbm$, 
\ie, that for every $M<\infty$ there exists a compact subset
$K\subset\cbbm$ such that 
\begin{equation}
  \label{km}
\limsup_{\eps\to0}\Bigpar{\eps\ln\P(\eps^{1/2}B\notin K)}\le -M;
\end{equation}
this exponential tightness is easily verified by choosing a $\gam$
with $\gb<\gam<1/2$ and taking
$K=\set{f\in\cbbm: \norm{f}_{C^\gam} \le L}$ for a large $L$.
We omit the verifications that $K$ is compact and satisfies \eqref{km}
if $L=L(M)$ is large enough.

The rest of proof of \refT{T1b} is entirely the same as for \refT{T1}.
\end{proof}

\begin{acks}
We thank Philippe Chassaing, Philippe Flajolet, and Nevin Kapur for
helpful comments.   
\end{acks}

\newcommand\AAP{\emph{Adv. Appl. Probab.} }
\newcommand\JAP{\emph{J. Appl. Probab.} }
\newcommand\JAMS{\emph{J. \AMS} }
\newcommand\MAMS{\emph{Memoirs \AMS} }
\newcommand\PAMS{\emph{Proc. \AMS} }
\newcommand\TAMS{\emph{Trans. \AMS} }
\newcommand\AnnMS{\emph{Ann. Math. Statist.} }
\newcommand\AnnPr{\emph{Ann. Probab.} }
\newcommand\CPC{\emph{Combin. Probab. Comput.} }
\newcommand\JMAA{\emph{J. Math. Anal. Appl.} }
\newcommand\RSA{\emph{Random Struct. Alg.} }
\newcommand\ZW{\emph{Z. Wahrsch. Verw. Gebiete} }
\newcommand\DMTCS{\jour{Discr. Math. Theor. Comput. Sci.} }

\newcommand\AMS{Amer. Math. Soc.}
\newcommand\Springer{Springer}
\newcommand\Wiley{Wiley}

\newcommand\vol{\textbf}
\newcommand\jour{\emph}
\newcommand\book{\emph}
\newcommand\inbook{\emph}
\def\no#1#2,{\unskip#2, no.~#1,} 

\newcommand\webcite[1]{\hfil\penalty0\texttt{\def~{\~{}}#1}\hfill\hfill}
\newcommand\webcitesvante{\webcite{http://www.math.uu.se/\~{}svante/papers/}}
\newcommand\arxiv[1]{\webcite{arXiv:#1}}

\def\nobibitem#1\par{}


\end{document}